\documentclass[a4paper,12pt]{article}

\usepackage{color}
\usepackage[dvips]{graphicx}
\usepackage{amscd}
\usepackage{amsmath}
\usepackage{ascmac}
\usepackage{cases}
\usepackage{amssymb}
\usepackage{amsfonts}
\usepackage{amsthm}
\usepackage[all]{xy}
\usepackage{enumerate}
\usepackage{wrapfig}

\theoremstyle{definition}
\newtheorem{Thm}{{\bf Theorem}}[section]

\newtheorem{Def}[Thm]{{\bf Definition}}
\newtheorem{Rem}[Thm]{{\bf Remark}}

\numberwithin{equation}{section}

\title{A nonclassical algebraic solution of 
a $3$-variable irregular Garnier system}
\author{Arata Komyo} 
\date{}

\begin{document}

\maketitle

\begin{abstract}
In this paper, a nonclassical algebraic solution of 
a $3$-variable irregular Garnier system is constructed. 
Diarra--Loray have studied classification of algebraic solutions of irregular Garnier systems.
There are two type of the algebraic solutions: classical type and pull-back type.
They have shown that there are exactly three nonclassical algebraic solutions 
for $N$-variables irregular Garnier systems with $N>1$.
Explicit forms of two of the three solutions are already given.
The solution constructed in the present paper is the remaining algebraic solution.



\end{abstract}

\section{Introduction}

The $N$-variable \textit{Garnier system} $G(1, \ldots ,1;N)$ 
is the completely integrable Hamiltonian system
$$
G(\underbrace{1, \ldots ,1}_{N+3};N) \colon
\left\{
\begin{aligned}
\frac{\partial q_j}{\partial t_i} &= \frac{\partial 
H_{t_i}(\boldsymbol{t},\boldsymbol{q},\boldsymbol{\eta})}{\partial \eta_j} \\
\frac{\partial \eta_j}{\partial t_i} &= -\frac{\partial 
H_{t_i}(\boldsymbol{t},\boldsymbol{q},\boldsymbol{\eta})}{\partial q_j} ,
\end{aligned}
\right.
$$
$i=0,1,2,\ldots,N-1$ and $j=1,2,3,\ldots,N$, where
\begin{equation*}
H_{t_i}(\boldsymbol{t},\boldsymbol{q},\boldsymbol{\eta})= -\frac{\Lambda(t_i)}{T'(t_i)}
\left[ \sum^{N}_{k=1} \frac{T(q_k)}{(q_k-t_i) \Lambda'(\nu_k)}
\left\{ \eta_k^2 -\sum^{N+1}_{m=0} \frac{\theta_m-\delta_{im}}{q_k - t_m} \eta_k
+ \frac{\kappa}{q_k(q_k -1)} \right\}
\right]
\end{equation*}
with $t_{N}:=0$, $t_{N+1}:=1$,
\begin{itemize}
\item $\theta_{0},\ldots, \theta_{N+1}, \theta_{\infty}$ are constant,
\item
$\kappa:= \frac{1}{4}\left\{ (\sum^{N+1}_{m=0} \theta_m -1)^2 - (\theta_{\infty}^2 +1) \right\}$,
\item $\Lambda(x):= \prod^{N}_{k=1} (x- q_k)$,  
$T(x):= \prod^{N+1}_{m=0} (x-t_m)$, 
and 
\item $\delta_{im}$ is the Kronecker delta
\end{itemize}
(see \cite{Garnier1}, \cite{Garnier2}, and \cite{Oka}).
The Garnier system is derived by the isomonodromic deformation of
a second order linear differential equation of the Fuchsian type with $N+3$ regular singular points,
$t_0,\ldots,t_{N-1}, 0,1,\infty$.
If $N=1$, the system $G(1, 1,1,1;1)$ is just $P_{VI}$, which is the Painlev\'e VI equation.

There are some methods to derive the Painelev\'e equations $P_J$ ($J=I,\ldots,V$).
For example, first,
$P_J$ ($J=I,\ldots,V$) are obtained from $P_{VI}$
by means of successive process of degeneration. 
Second, each of the Painelev\'e equations $P_J$ ($J=I,\ldots,V$)
is derived by 
the isomonodromic deformation of the degenerate linear differential equation.
This degenerate linear differential equation is given by 
the confluence of singularities of the linear differential equation.
Here the process of degeneration in the first method is
parallel with the confluence of singularities in the second method.
We can also consider the successive process of degeneration
for the Garnier system $G(1, \ldots ,1;N)$. 
For example, Kimura \cite{Kim} gave degenerations of $G(1, \ldots ,1;N)$ when $N=2$ 
(see also \cite{HKNS, Kawa}).
This process of degeneration is also
parallel with the confluence of singularities of the linear differential equation.
That is, the degenerations of $G(1,1,1,1 ,1;2)$ are the isomonodromy preserving deformation equations
of the degenerate linear differential equations
given by the confluence of singularities. 
Remark that the degenerate linear differential equations
have irregular singular points.
Suzuki considered the degeneration $G(1,\ldots,1,2;N)$ of $G(1, \ldots ,1;N)$,
which is the isomonodromy preserving deformation equations of linear differential equations
with $N+1$ regular singular points and an irreuglar singular point of Poincar\'e rank $1$.
Kawamuko (\cite{Kawa-2} and \cite{Kawa-1}) considered $G(1, N+2 ;N)$,
which is the isomonodromy preserving deformation equations of linear differential equations
with a regular singular point at $0$ and an irreuglar singular point at $\infty$ of Poincar\'e rank $N+1$.
We will call the isomonodromy preserving deformation equations of second order linear differential equations
with irregular singular points {\it irregular Garnier system}.
In the present paper, we are concerned with
the 3-variable irregular Garnier system $G(2,2,2;3)$.
This is the isomonodromy preserving deformation equations of second order linear differential equations
with 3 irregular singular points with Poincar\'e rank 1.

\subsection{Background and Motivation}
The Garnier system $G(1, \ldots ,1;N)$ has parameters $\theta_{0},\ldots, \theta_{N+1}, \theta_{\infty}$.
These parameter are related to the local exponents at singular points of the 
second order linear differential equation of the Fuchsian type with $N+3$ regular singular points.
In general, solutions of $G(1, \ldots ,1;N)$ are expected to be very transcendental.
However for $G(1, \ldots ,1;N)$ with special parameters,
there arise algebraic solutions.
When $N=1$, that is, $G(1, \ldots ,1;N)$ is the Painlev\'e equation $P_{VI}$,
Boalch \cite{Boa, Boa1,Boa2,Boa3, Boa4} has found the exhaustive list of algebraic solutions of $P_{VI}$.
The fact that this list is complete has been proved by Lisovyy--Tykhyy \cite{LT}.
Algebraic solutions of the Garnier system with $N>1$ have been investigated by several authors 
\cite{Dia, Cou, CoMo, CH, Gir, CaMa, Kom-1, DGS}.
But we are far from a complete classification of algebraic solutions of $G(1, \ldots ,1;N)$.

On the other hand, for irregular Garnier systems, 
there is a complete classification of algebraic solutions due to Diarra--Loray \cite{DL1}.
They have shown that any algebraic solution of an irregular 
Garnier system is one of the following types:
\begin{itemize}
\item classical type, which comes from the deformation of a system with diagonal or dihedral differential 
Galois group

\item pull-back type, which comes from the deformation obtained by 
pulling back a fixed system by an algebraic family of ramified covers.

\end{itemize}
(See \cite[Theorem 1 and Corollary 12]{DL1}).
This result is a direct consequence of
the classification theorem of flat meromorphic $\mathfrak{sl}_2$-connections on a 
projective manifold in \cite{LPT}.
Moreover, Diarra--Loray have studied a classification of algebraic solutions of pull-back type.
They have shown that $N$-variable irregular Garnier systems with $N>3$
admit only classical algebraic solutions.
Moreover, up to canonical transformations (see \cite{DM}),
there are exactly three nonclassical algebraic solutions for
$N$-variable irregular Garnier systems with $N>1$ (see \cite[Theorem 2]{DL1}).
Remark that these algebraic solutions are pull-back type.
For two of the three algebraic solutions, $N=2$.
When $N=2$, the list of explicit irregular Garnier systems is provided in \cite{Kim, Kawa}.
By finding fixed systems on $\mathbb{P}^1$ and  
algebraic families of ramified covers $\mathbb{P}^1 \rightarrow \mathbb{P}^1$,
Diarra--Loray \cite{DL1} have given explicit forms of those two algebraic solutions 
under Kimura's notation \cite{Kim}.
On the other hand,
for the last of the three algebraic solutions, $N=3$.
Diarra--Loray have found a fixed systems on $\mathbb{P}^1$ and  
an algebraic family of ramified covers $\mathbb{P}^1 \rightarrow \mathbb{P}^1$,
which correspond to the remaining algebraic solution. 
However, the corresponding explicit irregular Garnier system
had not been known.

On the other hand, there is a study of explicit description of isomonodromic deformations 
of second order linear differential equation with irregular singularities 
by using apparent singularities in \cite{Kom}.
So the purpose of this paper is 
to give the explicit irregular Garnier system corresponding to the remaining algebraic solution
by using the study in \cite{Kom}.
In consequence, we will give an explicit form of the algebraic solution
by using the fixed system on $\mathbb{P}^1$ and  
the algebraic family of ramified covers $\mathbb{P}^1 \rightarrow \mathbb{P}^1$
which are found by Diarra--Loray \cite[the first line of Table 2]{DL1}.

\subsection{Main Result}

We set $N=3$.
We define the Hamiltonians $H_{t_i}(\boldsymbol{t},\boldsymbol{q},\boldsymbol{\eta})$ $(i=0,1,2)$
as in \eqref{Hamiltonians} below.
In Section \ref{SubSectGarnier}, we will show that 
\begin{equation}\label{GarnierInIntro}
\frac{\partial q_j}{\partial t_i} =  \frac{\partial H_{t_i}(\boldsymbol{t},\boldsymbol{q},\boldsymbol{\eta})}{\partial \eta_j} ,
\qquad 
\frac{\partial \eta_j}{\partial t_i} = - \frac{\partial H_{t_i}(\boldsymbol{t},\boldsymbol{q},\boldsymbol{\eta})}{\partial q_j}
\end{equation}
($i=0,1,2$ and $j=1,2,3$)
is an irregular Garnier system.
More precisely, the system \eqref{GarnierInIntro}
is the isomonodromy preserving deformation equations of second order linear differential equations
with 3 irregular singular points with Poincar\'e rank 1.
Here we impose that these linear differential equations have a specific local normal form
in the sense of Hukuhara--Levelt--Turrittin for each singular point (as in \eqref{LocalFormalDataMatrix}).
Our main result is the following theorem.

\begin{Thm}\label{MainTheorem}
{\it We define $q_j ,\eta_j$ $(j=1,2,3)$ implicitly as follows}:
\begin{equation}\label{AlgebraicSol}
\left\{
\begin{aligned}
&q_j^3 +\frac{s_0-s_1-3\, s_2}{2\, s_2} q_j^2 
+\frac{-3\, s_0-s_1+ s_2}{2\, s_2}q_j +\frac{s_0}{s_2} =0 & j=1,2,3\\ 
&\eta_j=\frac{1}{3\, q_j} +\frac{1}{3\, (q_j-1)} & j=1,2,3 \\
&t_{i}^2  = s_i^3  & i=0,1,2.
\end{aligned}
\right.
\end{equation}
{\it Then such $q_j ,\eta_j$ $(j=1,2,3)$ satisfy the irregular Garnier system \eqref{GarnierInIntro}.}
\end{Thm}

The organization of this paper is as follows.
In Section \ref{Sect2}, we will derive the irregular Garnier system \eqref{GarnierInIntro}.
For this purpose, first, we will describe the corresponding second order linear differential equations.
Remark that we will consider second order linear differential equations 
as rank 2 connections on $\mathbb{P}^1$.
We will describe a family of rank 2 connections on $\mathbb{P}^1$
and consider the isomonodromic deformation of the rank 2 connections on $\mathbb{P}^1$.
In Section \ref{Sect3}, we will give an algebraic solution of the irregular Garnier system \eqref{GarnierInIntro}.
This is just \eqref{AlgebraicSol}.
This algebraic solution is constructed by the pull-back of
the fixed system on $\mathbb{P}^1$ by
the algebraic family of ramified covers $\mathbb{P}^1 \rightarrow \mathbb{P}^1$,
which are found by Diarra--Loray.
This system and the algebraic family will be described. 
Moreover, in Section \ref{2022_9_17_19_52}, 
we will give explicit description of 
the corresponding $\tau$-function,
which is defined as in \cite[Section 5]{JMU}.
In the appendix, we will try to recheck that the functions $q_j ,\eta_j$ defined by \eqref{AlgebraicSol} satisfy 
the irregular Garnier system \eqref{GarnierInIntro} by explicit calculation.

\section{The irregular Garnier system}\label{Sect2}

In this paper, we will treat only a specific irregular Garnier system.
Irregular Garnier systems come from the isomonodromic deformations
of rank 2 connections on $\mathbb{P}^1$ with irregular singular points. 
First we will describe rank 2 connections on $\mathbb{P}^1$ 
corresponding to our irregular Garnier system.
Let $x$ be a coordinate on an affine chart of $\mathbb{P}^1$.
We denote by $\infty$ the delated point and 
set $w:=1/x$.
Let $[0]$, $[1]$, and $[\infty]$ be divisors on $\mathbb{P}^1$
defined by $x=0,x=1$, and $w=0$, respectively. 
We set $D:= 2[0] + 2[1] +2[\infty]$.
We denote by $E_k$ the rank 2 vector bundle
$\mathcal{O}\oplus \mathcal{O} (k)$ on $\mathbb{P}^1$ for an integer $k$.
Let $\nabla$ be a connection on $E_1$ whose polar divisor is $D$:
\begin{equation}\label{OrigConn}
\nabla \colon E_1 \longrightarrow E_1 \otimes \Omega^1_{\mathbb{P}^1} (D).
\end{equation}
The restrictions $\nabla|_{2[0]}$, $\nabla|_{2[1]}$, $\nabla|_{2[\infty]}$
at the divisors $2[0]$, $2[1]$, $2[\infty]$ 
are well-defined, respectively.
We assume that the points $0$, $1$, $\infty$ are
unramified irregular singular points of $\nabla$.
So we can diagonalize
$\nabla|_{2[0]}$, $\nabla|_{2[1]}$, $\nabla|_{2[\infty]}$
by some local gauge transformations at $0$, $1$, $\infty$, respectively.
We call the diagonal entries of the diagonalizations for the points $0,1,\infty$
{\it local formal data of} $\nabla$.
We are interested in an algebraic solution of an irregular Garnier system.
To have an algebraic solution,
we impose that the connections \eqref{OrigConn}
have the following special local formal data: 
\begin{equation}\label{LocalFormalData}
\left(
\left( \epsilon \frac{2\, t_0 dx}{x^2} -\frac{dx}{6\, x}   \right)_{\epsilon=+,-} ,
\left( \epsilon \frac{2\, t_1dx}{(x-1)^2} -\frac{dx}{6\, (x-1)} \right)_{\epsilon=+,-} ,
\left( \epsilon \frac{2\, t_2dw}{w^2} -\frac{dw}{6\, w}   \right)_{\epsilon=+,-} 
\right) .
\end{equation}
Here $t_i$ is a complex number with $t_i\not=0$ for each $i=0,1,2$.
That is, 
$\nabla|_{2[0]}$, $\nabla|_{2[1]}$, $\nabla|_{2[\infty]}$
are diagonalized as 
\begin{equation}\label{LocalFormalDataMatrix}
\begin{aligned}
&
\begin{pmatrix}
2\, t_0& 0 \\
0 & -2\, t_0 
\end{pmatrix}\frac{dx}{x^2}
+\begin{pmatrix}
-\frac{1}{6}& 0 \\
0 & -\frac{1}{6}
\end{pmatrix}\frac{dx}{x} \\
&
\begin{pmatrix}
2\, t_1& 0 \\
0 & -2\, t_1 
\end{pmatrix}\frac{dx}{(x-1)^2}
+\begin{pmatrix}
-\frac{1}{6}& 0 \\
0 & -\frac{1}{6}
\end{pmatrix}\frac{dx}{x-1}\\
&
\begin{pmatrix}
2\, t_2& 0 \\
0 & -2\, t_2 
\end{pmatrix}\frac{dw}{w^2}
+\begin{pmatrix}
-\frac{1}{6}& 0 \\
0 & -\frac{1}{6}
\end{pmatrix}\frac{dw}{w},
\end{aligned}
\end{equation}
by some local gauge transformations.
Then the corresponding irregular Garnier system 
has an algebraic solution that we aim at.

\begin{Rem}
The connection \eqref{OrigConn} is not trace-free.
But we can transform this connection into a trace-free connection.
Indeed, we apply some birational bunlde transformation 
$E_1 \dashrightarrow E_0$
for the connection \eqref{OrigConn}, and
take the tensor product of this birational bunlde transformation and 
a rank 1
connection $\mathcal{O} \rightarrow \mathcal{O} \otimes 
\Omega^1_{\mathbb{P}^1}([0]+[1]+[\infty])$.
Then we may give a trace-free connection
with the local formal data 
\begin{equation*}
\left(
\left( \epsilon \frac{2\, t_0 dx}{x^2}    \right)_{\epsilon=+,-} ,
\left( \epsilon \frac{2\, t_1dx}{(x-1)^2}  \right)_{\epsilon=+,-} ,
\left( \epsilon \left( \frac{2\, t_2dw}{w^2} + \frac{dw}{2\, w} \right)   \right)_{\epsilon=+,-} 
\right) .
\end{equation*}
\end{Rem}

We will consider the isomonodromic deformations of the connections \eqref{OrigConn}
along the parameters $t_0$, $t_1$, and $t_2$.
So the time variables of our irregular Garnier system are $t_0$, $t_1$, and $t_2$.
The space of time variables is $(\mathbb{C}^*)^3$.
We consider moduli spaces of such connections (with some generic conditions).
When we fix the time variables (that is, local formal data), 
we call the moduli space {\it the moduli space with fixed time variables}.
When the time variables may vary, 
we call the moduli space {\it the extended moduli space} 
(for details, for example, see \cite[Section 2.4]{Kom}).
The isomonodromic deformations mean vector fields on 
the extended moduli space.
Remark that the dimension of our moduli space
with fixed time variables is $6$.
To drive an irregular Garnier system, we will introduce a coordinate system 
on the moduli space with fixed time variables and will describe 
the isomonodromic deformations by using this coordinate system.
We will use the coordinate system given by {\it apparent singularities}
(for details, see \cite[Section 1]{DL2} and \cite[Section 2]{Kom}).
By using the coordinate system, we have
a birational map from this moduli space
with fixed time variables to $\mathrm{Sym}^3(\mathbb{C}^2)$
(see \cite[Theorem 1.1]{DL2} and \cite[Section 2.4]{Kom}).
Here $\mathrm{Sym}^3(\mathbb{C}^2)$ is
the 3-fold symmetric product of $\mathbb{C}^2$. 
We denote
by $\{ (q_1,p_1),(q_2,p_2),(q_3,p_3) \}$ a point on $\mathrm{Sym}^3(\mathbb{C}^2)$. 
Now we describe this birational map
by introducing a family of the connections \eqref{OrigConn} parametrized 
by the Zariski open subset $\mathcal{M}$ of $\mathrm{Sym}^3(\mathbb{C}^2)$.
Here $\mathcal{M}$ is defined by
$$
\mathcal{M}= 
\left\{ \{ (q_1,p_1),(q_2,p_2),(q_3,p_3) \} \in \mathrm{Sym}^3(\mathbb{C}^2)
\ \middle| \ 
\begin{array}{l}
\text{$q_j \not\in \{0,1,\infty\}$ ($j=1,2,3$)} \\
\text{$q_{j_1}\not= q_{j_2}$ ($j_1 \not= j_2$)}
\end{array}
\right\}.
$$
This family will give an inverse map of the birational map.
That family has constructed in \cite{DL2}.
Now, we recall the construction of that family.
To construct that family, we use a family of connections
\begin{equation}\label{ConnNormal}
\nabla' \colon E_4 \longrightarrow E_4 \otimes \Omega^1_{\mathbb{P}^1} (D+q_1+q_2+q_3)
\end{equation}
parametrized by $\mathcal{M}$.
Here we assume that $q_1,q_2,q_3$ are apparent singular points of $\nabla'$,
that is, these singularities of $\nabla'$ 
can be removed by a birational bundle transformation.

Now we recall the family of connections $\nabla'$. 
We set $U_0 :=\mathbb{P}^1\setminus 0$ and 
$U_\infty :=\mathbb{P}^1\setminus \infty$.
We fix trivializations 
$$
(U_0,\ \varphi_{0}\colon E_4|_{U_0} \rightarrow \mathcal{O}_{U_0}^{\oplus 2}) 
\quad \text{and} \quad 
(U_\infty, \  \varphi_{\infty}\colon E_4|_{U_\infty} \rightarrow \mathcal{O}_{U_\infty}^{\oplus 2})
$$
of $E_4$ such that $\varphi_{0} \circ \varphi_{\infty}^{-1}=
\begin{pmatrix} 1 & 0 \\ 0 & x^4 \end{pmatrix}$.
As in \cite[Section 3]{DL2} and \cite[Section 2.2]{Kom},
we define a family $\Omega^{(4)}$ of connection matrices with respect to 
$(U_0,\ \varphi_{0}) $ parametrized by $\mathcal{M}$
as 
\begin{equation}\label{NormalForm}
\Omega^{(4)}=
\begin{pmatrix}
0& \frac{1}{x^2(x-1)^2} \\
c_2(x;\boldsymbol{q}, \boldsymbol{p}) & d_2(x;\boldsymbol{q}, \boldsymbol{p})
\end{pmatrix}\, dx,
\end{equation}
where $c_2$ and $d_2$ have the following forms:
\begin{equation*}
\begin{aligned}
c_2(x;\boldsymbol{q}, \boldsymbol{p})&:=
\frac{C_0(x)}{x^2}
+\frac{C_1(x)}{(x-1)^2}+\sum_{i=1}^3\frac{p_i}{x-q_i} +z^3C_{\infty}(x)\\
&\qquad 
+\tilde{C}_{q_1}(x-q_2)(x-q_3)
+\tilde{C}_{q_2}(x-q_1)(x-q_3)
+\tilde{C}_{q_3}(x-q_2)(x-q_1)  , \quad \text{ and} \\
d_2(x;\boldsymbol{q}, \boldsymbol{p})&:=
 \frac{D_0(x)}{x^2}
+\frac{D_1(x)}{(x-1)^2}
- \sum_{i=1}^2 \frac{1}{x-q_i} +D_{\infty}(x).
\end{aligned}
\end{equation*}
Here $C_0,C_1,C_\infty$, 
$D_0,D_1$, and $D_\infty$ are determined as follows:
\begin{equation}\label{PolyCD}
\left\{
\begin{aligned}
&C_0(x) =  4 \, t_0^2 (1-2\, x )  \\
&C_1(x) = 4\,  t_1^2 ( 2\, x-1) \\
&C_\infty(x) = 4 \,  t_2^2 (x -2) 
\end{aligned} 
\right.
\qquad 
\left\{
\begin{aligned}
&D_0(x) =  -\frac{x}{3} \\
&D_1(x) = -\frac{x-1}{3} \\
&D_\infty(x) = 0
\end{aligned}
\right..
\end{equation}
We can extend connections $d+ \Omega^{(4)} $ on $U_0$
to connections on $E_4$ over $\mathbb{P}^1$, naturally.
We denote by $\nabla_4^{\text{fam}}(\boldsymbol{q}, \boldsymbol{p})$ the extended connections.
By the condition \eqref{PolyCD}, the points $0,1,\infty$ on $\mathbb{P}^1$
are unramified irregular singular points of $\nabla_4^{\text{fam}}(\boldsymbol{q}, \boldsymbol{p})$. 
These local formal data are just \eqref{LocalFormalData}
(for detail, see \cite[Lemma 4.1]{DL2}).
Moreover we define the constants $\tilde{C}_{q_j}$ ($j=1,2,3$) in \eqref{NormalForm} 
as follows:
$$
{\small
\begin{aligned}
\tilde{C}_{q_j}= \frac{1}{Q'(q_j)}
\Biggl(  \frac{p_j^2}{q_j^2(q_j-1)^2} 
&+\frac{ q_j p_j + 12\,  t_0^2 (2\, q_j-1 )  }{ 3\,  q_j^2}
+\frac{(q_j-1) p_j - 12\,  t_1^2 ( 2\, q_j-1)  }{ 3\, (q_j-1)^2} \\
&\quad +\sum_{k\in\{1,2,3\} \setminus \{ j\}} \frac{p_j-p_k}{q_j-q_k} - 4\, t_2^2 q_j^3 (q_j -2)  
\Biggr),
\end{aligned}}
$$
where $Q'(x)=(x-q_1)(x-q_2)+(x-q_2)(x-q_3)+(x-q_3)(x-q_1)$ (see \cite[the formula (5.5)]{Kom}).
Then the poles $q_1,q_2,q_3$ of the connections 
$d+ \Omega^{(4)} $ are apparent singular points.
The connections $\nabla_4^{\text{fam}}(\boldsymbol{q}, \boldsymbol{p})$
are just the family of connections $\nabla'$ that we will use. 
Since we assume that $q_j \not\in \{0,1,\infty\}$ ($j=1,2,3$),
there exists a birational bundle transformation $ \phi \colon E_1 \dashrightarrow E_4$
such that the pull-back $ \phi^* \nabla_4^{\text{fam}}(\boldsymbol{q}, \boldsymbol{p})$
is a family of connections \eqref{OrigConn}
parametrized by $\mathcal{M}$
(see \cite[the formula (1.5)]{DL2}
and \cite[Proposition 2.4]{Kom}).
We have a map 
$$(\boldsymbol{q}, \boldsymbol{p})=
\{ (q_1,p_1),(q_2,p_2),(q_3,p_3) \}
\longmapsto 
 \phi^* \nabla_4^{\text{fam}}(\boldsymbol{q}, \boldsymbol{p})
$$
from $\mathcal{M}$ to the moduli space with fixed time variables.
This is the inverse map of our birational map 
from the moduli space with fixed time variables
to $\mathrm{Sym}^3(\mathbb{C}^2)$.

Remark that we may consider 
$ \phi^* \nabla_4^{\text{fam}}(\boldsymbol{q}, \boldsymbol{p})$
as a family of connections \eqref{OrigConn}
parametrized by $\mathcal{M} \times (\mathbb{C}^*)^3$.
Here $(\mathbb{C}^*)^3$ is the space of time variables.
So we can extend the birational map 
to a birational map between 
$\mathcal{M} \times (\mathbb{C}^*)^3$ and 
the extended moduli space.
The natural projection from the extended moduli space to $(\mathbb{C}^*)^3$
is compatible with the projection 
$\mathcal{M} \times (\mathbb{C}^*)^3 \rightarrow (\mathbb{C}^*)^3$
through that birational map.

\subsection{The corresponding Garnier system}\label{SubSectGarnier}

Now,
we consider the isomonodromic deformations of the connections \eqref{OrigConn}
along the parameters $t_0, t_1$, and $t_2$.
The parameter space of $t_0, t_1$, and $t_2$ is $(\mathbb{C}^*)^3$.
For vector fields $\frac{\partial}{\partial t_0}$, 
$\frac{\partial}{\partial t_1}$, $\frac{\partial}{\partial t_2}$ on $(\mathbb{C}^*)^3$, 
we can define lifts of these vector fields under the natural projection 
from the extended moduli space to $(\mathbb{C}^*)^3$.
Here the lifts 
mean the isomonodromic deformations of the connection corresponding to each point on
the extended moduli space.
By the birational map, 
we have the vector fields on $\mathcal{M}\times (\mathbb{C}^*)^3$
corresponding to the lifts.
We give explicit description of
these vector fields on $\mathcal{M}\times (\mathbb{C}^*)^3$.
They are just our irregular Garnier system.

To describe the vector fields of the isomonodromic deformations explicitly,
we will use $q_j$ and $\eta_j$ ($j=1,2,3$) as coordinates on $\mathcal{M}$.
Here $\eta_j$ ($j=1,2,3$) is defined by
\begin{equation}\label{ETA}
\left\{
\begin{aligned}
&\eta_1 =\frac{p_1}{q_1^2(q_1-1)^2} +\frac{1}{3\, q_1} +\frac{1}{3(q_1-1)} \\
&\eta_2 =\frac{p_2}{q_2^2(q_2-1)^2} +\frac{1}{3\, q_2} +\frac{1}{3(q_2-1)} \\
&\eta_3 =\frac{p_3}{q_3^2(q_3-1)^2} +\frac{1}{3\, q_3} +\frac{1}{3(q_3-1)}
\end{aligned}
\right..
\end{equation}

\begin{Rem}
The parameters $\eta_j$ ($j=1,2,3$) defined in \eqref{ETA} come from 
$$
\eta_j = \frac{p_j}{q_j^2(q_j-1)^2} - \frac{D_0(q_j)}{q_j^2} - \frac{D_1(q_j)}{(q_j-1)^2} -
D_\infty(q_j) \qquad (j=1,2,3),
$$
which are introduced in \cite[Corollary 3.15]{Kom}.
Now the coefficients of $D_0(x)$, $D_1(x)$, and $D_\infty(x)$
are independent of $t_0,t_1,t_2$. 
It means that the trace of the local formal data at each unramified  irregular singular point 
is independent of $t_0,t_1,t_2$.
Then we have the following equalities of $2$-forms on $\mathcal{M}\times (\mathbb{C}^*)^3$:
$$
\begin{aligned}
d \eta_j \wedge d q_j 
&= d\left( 
\frac{p_j}{q_j^2(q_j-1)^2} - \frac{D_0(q_j)}{q_j^2(q_j-1)^2} - \frac{D_1(q_j)}{q_j^2(q_j-1)^2} -
D_\infty(q_j)
\right) \wedge d q_j  \\
&= d\left( 
\frac{p_j}{q_j^2(q_j-1)^2}  \right) \wedge d q_j .
\end{aligned}
$$
So we can use $\frac{p_j}{q_j^2(q_j-1)^2}$ ($j=1,2,3$) as parameters to describe the Hamilton system 
of the irregular Garnier system instead of \eqref{ETA}.
\end{Rem}

Now we will define Hamiltonians of our irregular Garnier system.
We fix compatible framings $\Phi_0$, $\Phi_1$, $\Phi_\infty$
 at $0$, $1$, $\infty$, respectively, as follows:
$$
\Phi_0 = 
\begin{pmatrix}
\frac{1}{2\, t_1} & - \frac{1}{2\, t_1} \\
1 & 1 
\end{pmatrix}, \qquad 
\Phi_1 = 
\begin{pmatrix}
\frac{1}{2\, t_2} & - \frac{1}{2\, t_2} \\
1 & 1 
\end{pmatrix}, \qquad 
\Phi_{\infty} = 
\begin{pmatrix}
-\frac{1}{2\, t_3} & \frac{1}{2\, t_3} \\
1 & 1 
\end{pmatrix}.
$$
That is, by these matrices, we can diagonalize the leading coefficients of the Laurent expansions of 
\eqref{NormalForm} at 
the unramified irregular singular points $0$, $1$, and $\infty$.
For each unramified irregular singular point, 
there exists a unique formal power series 
$\mathrm{id}+ \sum_{k=1}^{\infty} \Xi_k \tilde{x}^k$
such that 
\begin{itemize}
\item For each $k=1,2,3,\ldots$,
the coefficient $\Xi_k$ is a two-by-two matrix
and the diagonal entries of $\Xi_k$ are zero, and
\item the Laurent expansion of 
\eqref{NormalForm} at the singular point is diagonalized by
$\Phi(\mathrm{id}+ \sum_{k=1}^{\infty} \Xi_k \tilde{x}^k)$,
where $\Phi$ is the fixed compatible framing at the point.
\end{itemize}
Here we set $\tilde{x}:=x$, $\tilde{x}:=(x-1)$, and $\tilde{x}:=w$ 
when the unramified irregular singular point which 
we consider is $0$, $1$, and $\infty$, respectively.
We consider these diagonalizations until the $\tilde{x}^0$-terms: 
\begin{equation*}
\begin{aligned}
&\begin{pmatrix}
2\, t_0& 0 \\
0 & -2\, t_0 
\end{pmatrix}\frac{dx}{x^2}
+\begin{pmatrix}
-\frac{1}{6}& 0 \\
0 & -\frac{1}{6}
\end{pmatrix}\frac{dx}{x}
+\begin{pmatrix}
\theta_{0}^+(\boldsymbol{t},\boldsymbol{q},\boldsymbol{\eta})& 0 \\
0 &\theta_{0}^-(\boldsymbol{t},\boldsymbol{q},\boldsymbol{\eta})
\end{pmatrix}dx+O(x) \\
&\begin{pmatrix}
2\, t_1& 0 \\
0 & -2\, t_1 
\end{pmatrix}\frac{dx}{(x-1)^2}
+\begin{pmatrix}
-\frac{1}{6}& 0 \\
0 & -\frac{1}{6}
\end{pmatrix}\frac{dx}{x-1}
+\begin{pmatrix}
\theta_{1}^+(\boldsymbol{t},\boldsymbol{q},\boldsymbol{\eta})& 0 \\
0 &\theta_{1}^-(\boldsymbol{t},\boldsymbol{q},\boldsymbol{\eta})
\end{pmatrix}dx+O(x-1)\\
&\begin{pmatrix}
2\, t_2& 0 \\
0 & -2\,  t_2 
\end{pmatrix}\frac{dw}{w^2}
+\begin{pmatrix}
-\frac{1}{6}& 0 \\
0 & -\frac{1}{6}
\end{pmatrix}\frac{dw}{w}
+\begin{pmatrix}
\theta_{\infty}^+(\boldsymbol{t},\boldsymbol{q},\boldsymbol{\eta})& 0 \\
0 &\theta_{\infty}^-(\boldsymbol{t},\boldsymbol{q},\boldsymbol{\eta})
\end{pmatrix}dw+O(w).
\end{aligned}
\end{equation*}
Here we may compute 
$\theta_{0}^\pm(\boldsymbol{t},\boldsymbol{q},\boldsymbol{\eta})$,
$\theta_{1}^\pm(\boldsymbol{t},\boldsymbol{q},\boldsymbol{\eta})$, and
$\theta_{\infty}^\pm(\boldsymbol{t},\boldsymbol{q},\boldsymbol{\eta})$,
explicitly.
But we omit to describe these explicit forms.
We define Hamiltonians 
$H_{t_0} (\boldsymbol{t},\boldsymbol{q},\boldsymbol{\eta})$,
$H_{t_1} (\boldsymbol{t},\boldsymbol{q},\boldsymbol{\eta})$, and
$H_{t_2} (\boldsymbol{t},\boldsymbol{q},\boldsymbol{\eta})$ as follows.
\begin{Def}
We set 
$$
\left\{
\begin{aligned}
&H_{t_0} (\boldsymbol{t},\boldsymbol{q},\boldsymbol{\eta}):=
2\,  \theta_0^-(\boldsymbol{t},\boldsymbol{q},\boldsymbol{\eta})
-2\, \theta_0^+(\boldsymbol{t},\boldsymbol{q},\boldsymbol{\eta}), \\
&H_{t_1} (\boldsymbol{t},\boldsymbol{q},\boldsymbol{\eta}):=
2\,  \theta_1^-(\boldsymbol{t},\boldsymbol{q},\boldsymbol{\eta})
-2\,  \theta_1^+(\boldsymbol{t},\boldsymbol{q},\boldsymbol{\eta}), \\
&H_{t_2} (\boldsymbol{t}, \boldsymbol{q},\boldsymbol{\eta}):=
2\,  \theta_\infty^-(\boldsymbol{t},\boldsymbol{q},\boldsymbol{\eta})
-2\,  \theta_\infty^+(\boldsymbol{t},\boldsymbol{q},\boldsymbol{\eta}).
\end{aligned}
\right.
$$
\end{Def}
By explicit computation using
some computing environment (for example {\it Maple}), we have the following description of the Hamiltonians:
\begin{equation}\label{Hamiltonians}
{\footnotesize
\left\{
\begin{aligned}
&H_{t_0} (\boldsymbol{t},\boldsymbol{q},\boldsymbol{\eta})\\
&=- \frac{{q_{{1}}} q_{{2}}q_{{3}}}{ t_{{1}} } \sum_{j=1}^3 \left(
{\frac { q_{{j}} \left( q_{{j}}-1 \right) ^{2}}{ Q'(q_j)}  }{{ \eta_j}}^{2}
-{\frac {\left( 5\,{q_{{j}}}^{2}-9\,q_{{j}}+4 \right)}{3 \, Q'(q_j)}}{ \eta_j}
\right) 
+\sum_{j=1}^3{\frac {4\, t_{{1}}}{{q_{{j}}}^{2}}}
+{\frac {144\,{t_{{1}}}^{2}+144\,{t_{{2}}}^{2}-13}{36\,t_{{1}}}}  \\
&\quad 
+{\frac {4\, {t_{{3}}}^{2} \sigma_{{3}} \left( \sigma_{{1}}-2 \right) }{t_{{1}}}}
-{\frac {4\, t_{{1}} \left( 2\,\sigma_{{2}}-\sigma_{{1}} \right) }{q_{{1}}q_{{2}}q_{{3}}}}
-{\frac { 4\,{t_{{2}}}^{2} \left( {\sigma_{{1}}}^{2}-2\,\sigma_{{1}}\sigma_{{2}
}+3\,\sigma_{{1}}\sigma_{{3}}+{\sigma_{{2}}}^{2}-2\,\sigma_{{2}}\sigma
_{{3}}-2\,\sigma_{{1}}+2\,\sigma_{{2}}-4\,\sigma_{{3}}+1 \right) 
 }{ t_{{1}} \left( q_{{1}}-1 \right) ^{2} \left( q_{{2}}-1 \right) ^{2} \left( q_{{3}}-1 \right) ^{2}}} \\
&H_{t_1} (\boldsymbol{t},\boldsymbol{q},\boldsymbol{\eta}) \\
&=
-\frac{ \left( q_{{1}}-1 \right)  \left( q_{{2}}-1 \right)  \left( q_{{3}}-1 \right)}{t_{{2}}} \sum_{j=1}^3
\left(
{\frac {{q_{{j}}}^{2} \left( q_{{j}}-1 \right)   }{ Q'(q_j)}}{{ \eta_j}}^{2}
-\frac { q_{{j}} \left( 5 \, q_{{j}}-1 \right)  }{ 3\, Q'(q_j)  }{\eta_j}
\right) 
+\sum_{j=1}^3 {\frac {4\, t_{{2}}}{ \left( q_{{j}}-1 \right) ^{2}}}
+{\frac {144\,{t_{{1}}}^{2}+144\,{t_{{2}}}^{2}-13}{36\,t_{{2}}}} \\
&\quad 
+{\frac {4\, {t_{{3}}}^{2} (\sigma_{{1}}-\sigma_{{2}}+\sigma_{{3}}-1) \left( \sigma_{{1}}-1 \right) }{t_{{2}}}}
+{\frac {4\,t_{{2}} \left( 2\,\sigma_{{2}} -3\,\sigma_{{1}}+3 \right) }{ \left( q_{{1}}-1 \right)  \left( q_{{2}}-1 \right)  \left( q_{{3}}-1 \right) }}
-{\frac {4\,{t_{{1}}}^{2} \left( \sigma_{{1}}\sigma_{{2}}-\sigma_{{1}}\sigma_{{
3}}-{\sigma_{{2}}}^{2}+2\,\sigma_{{2}}\sigma_{{3}}-\sigma_{{2}}+\sigma
_{{3}} \right)  }{t_{{2}}{q_{{1}}}^{2}{q_{{2}}}^{2}{q_{{3}}}^{2}}}  \\
&H_{t_2} (\boldsymbol{t},\boldsymbol{q},\boldsymbol{\eta})\\
&=
 - \frac{1}{t_3} \sum_{j=1}^3
\left( 
{\frac {{q_{{j}}}^{2} \left( q_{{j}}-1 \right) ^{2}}{Q'(q_j)}}{{\eta_j}}^{2}
-{\frac {q_{{j}} \left( 2\,{q_{{j}}}^{2}-3\,q_{{j}}+1 \right)}{ 3\, Q'(q_j)  }} {\eta_j} \right)
+4\,{t_{{3}}} ( {\sigma_{{1}}}^2 -2\,\sigma_{{1}} - \sigma_{{2}}+1) 
-\frac {1}{36\,t_{{3}}} \\
&\quad 
-{\frac {4\,{t_{{1}}}^{2} \left( 2\,\sigma_{{3}}-\sigma_{{2}} \right) }{t_{{3}} {q_{{1}}}^{2}{q_{{2}}}^{2}{q_{{3}}}^{2}}} 
+{\frac {4\,{t_{{2}}}^{2} \left( 2\,\sigma_{{3}}-\sigma_{{2}}+1 \right) }{ t_{{3}} \left( q_{{1}}-1 \right) ^{2} \left( q_{{2}}-1 \right) ^{2} \left( q_{{3}}-1 \right) ^{2}}}
\end{aligned}
\right.,
}
\end{equation}
where $Q'(x)=(x-q_1)(x-q_2)+(x-q_2)(x-q_3)+(x-q_3)(x-q_1)$.
Here we used the elementary symmetric polynomials
$\sigma_1 =q_1+q_2+q_3$, 
$\sigma_2 =q_1q_2+q_2q_3+q_3q_1$, and
$\sigma_3 =q_1q_2q_3$
in the numerators of some terms in \eqref{Hamiltonians}.

Let $\hat\omega$ be the isomonodromy 2-form on 
the extended moduli space.
Remark that 
the vector fields of isomonodromic deformations are defined on 
the extended moduli space.
The isomonodromy 2-form means that $\hat\omega$ is a $2$-form on the extended moduli space
such that 
\begin{itemize}
\item the restriction of 
$\hat\omega$ to each fiber of the projection 
from the extended moduli space to the space of time variables is a symplectic form, 
and 
\item the interior products of $\hat\omega$
with the vector fields of isomonodromic deformations with respect to 
$\frac{ \partial}{\partial t_0}$ $\frac{ \partial}{\partial t_1}$, and $\frac{ \partial}{\partial t_2}$
vanish.
\end{itemize}
Now we have the birational map 
from $\mathcal{M} \times (\mathbb{C}^*)^3$
to the extended moduli space.
This birational map is described by the 
family of connections $ \phi^* \nabla_4^{\text{fam}}(\boldsymbol{q}, \boldsymbol{p})$ 
parametrized by $\mathcal{M} \times (\mathbb{C}^*)^3$.
We also denote by $\hat\omega$ the pull-back of the isomonodromy 2-form
by the birational map.

We will describe
the 2-form $\hat\omega$ on $\mathcal{M} \times (\mathbb{C}^*)^3$ explicitly.
By \cite[Theorem 3.14]{Kom},
there exist functions 
$f_{i_1,i_2} (\boldsymbol{t},\boldsymbol{q},\boldsymbol{\eta})$ 
on $\mathcal{M} \times (\mathbb{C}^*)^3$ such that 
$\hat\omega$ is described as follows:
$$
\begin{aligned}
\hat\omega &= \sum_{j=1}^3  d\eta_j \wedge dq_j 
- d \theta_0^+ \wedge d(2t_0) 
- d \theta_1^+ \wedge d(2t_1)
- d \theta_\infty^+ \wedge d(2t_2)  \\
&\qquad 
- d \theta_0^- \wedge d(-2t_0)
-d \theta_1^- \wedge d(-2t_1)
-d \theta_\infty^- \wedge d(-2t_2) 
+\sum_{i_1 <i_2} 
f_{i_1,i_2} (\boldsymbol{t},\boldsymbol{q},\boldsymbol{\eta}) \, dt_{i_1} \wedge dt_{i_2} .
\end{aligned}
$$
We may describe the isomonodromy 2-form by using the Hamiltonian \eqref{Hamiltonians}:
$$
\begin{aligned}
\hat\omega &= \sum_{j=1}^3  d\eta_j \wedge dq_j  
+ d H_{t_0} (\boldsymbol{t},\boldsymbol{q},\boldsymbol{\eta}) \wedge dt_0
+ d H_{t_1} (\boldsymbol{t},\boldsymbol{q},\boldsymbol{\eta}) \wedge dt_1
+ d H_{t_2} (\boldsymbol{t},\boldsymbol{q},\boldsymbol{\eta}) \wedge dt_2 \\
&\qquad +\sum_{i_1 <i_2} 
f_{i_1,i_2} (\boldsymbol{t},\boldsymbol{q},\boldsymbol{\eta}) \, dt_{i_1} \wedge dt_{i_2}.
\end{aligned}
$$
By the explicit form of $\hat\omega$, we have the following 
explicit form of the vector field of isomonodromic deformation with respect to 
$\frac{ \partial}{\partial t_i}$ for each $i=0,1,2$:
$$
\frac{ \partial}{\partial t_i} + \sum_{j=1}^3\left(
-\frac{\partial H_{t_i} (\boldsymbol{t},\boldsymbol{q},\boldsymbol{\eta})}{\partial \eta_j}
\frac{ \partial}{\partial q_j}
+\frac{\partial H_{t_i} (\boldsymbol{t},\boldsymbol{q},\boldsymbol{\eta})}{\partial q_j}
\frac{ \partial}{\partial \eta_j} \right)
$$
(for detail, see \cite[Corollary 3.15]{Kom}).
Finally, we have an irregular Garnier system by this explicit form.
\begin{Def}\label{Garnier_Def}
Let $H_{t_i}(\boldsymbol{t},\boldsymbol{q},\boldsymbol{\eta})$
be the Hamiltonian as in \eqref{Hamiltonians}.
We say the following system
\begin{equation}\label{GarnierInSect2}
\left\{
\begin{aligned}
\frac{\partial q_j}{\partial t_i} &=  \frac{\partial H_{t_i}(\boldsymbol{t},\boldsymbol{q},\boldsymbol{\eta})}{\partial \eta_j} \\
\frac{\partial \eta_j}{\partial t_i} &= - \frac{\partial H_{t_i}(\boldsymbol{t},\boldsymbol{q},\boldsymbol{\eta})}{\partial q_j}
\end{aligned}
\right. \qquad \text{for any $i=0,1,2$, $j=1,2,3$}
\end{equation}
the {\it irregular Garnier system with respect to the local formal data \eqref{LocalFormalData}}.
\end{Def}

By computation using 
the explicit formula \eqref{Hamiltonians} and
using some computing environment (for example {\it Maple}),
we have the following equality:
\begin{equation}\label{2022_9_16_23_08}
\frac{\partial H_{t_i}}{\partial t_j}
-\frac{\partial H_{t_j}}{\partial t_i}
-\{ H_{t_i} , H_{t_j} \}  =0
\end{equation}
for $i ,j= 0,1,2$ (where $i\neq j$).
Here $\{ H_{t_i} , H_{t_j} \} $ is a Poisson bracket of $H_{t_i}$ and $H_{t_j}$, 
that is, 
$$
\{ H_{t_i} , H_{t_j} \} 
= \sum_{k=1}^{3} 
\left(\frac{\partial H_{t_i}}{\partial q_k}\frac{\partial H_{t_j}}{\partial \eta_k}
-\frac{\partial H_{t_j}}{\partial q_k}\frac{\partial H_{t_i}}{\partial \eta_k}\right).
$$
The equality \eqref{2022_9_16_23_08}
will be used in Section \ref{2022_9_17_19_52}.

\section{Isomonodromic family and algebraic solution}\label{Sect3}

Let $U$ be an open subset of $(\mathbb{C}^*)^3$.
In this section, we will construct an isomonodromic family of connections \eqref{OrigConn}
parametrized by $U$ to give a solution of the irregular Garnier system \eqref{GarnierInSect2}.
That is, we construct a map
\begin{equation}\label{SolutionGeneral}
\begin{aligned}
l\colon  U &\longrightarrow  
\mathcal{M} \times (\mathbb{C}^*)^3 \\
\tilde {\boldsymbol{t}}=(\tilde t_0,\tilde t_1,\tilde t_2) &\longmapsto
 ( \{ (q_j(\tilde{\boldsymbol{t}}) , \eta_j(\tilde{\boldsymbol{t}})) \}_{j=1,2,3},
 \tilde {\boldsymbol{t}})
\end{aligned}
\end{equation}
such that
\begin{itemize}
\item the local formal data of the pull-back
$(l\times \mathrm{id})^* ( E_1 ,\, \phi^* \nabla_4^{\text{fam}}(\boldsymbol{q} , \boldsymbol{p}))$
are \eqref{LocalFormalData} with $t_0=\tilde  t_0, t_1=\tilde  t_1, t_2=\tilde  t_2$,

\item the family of equivalence classes of the (generalized)
monodromy representations of
the connections
$(l\times \mathrm{id})^* ( E_1 ,\, \phi^* \nabla_4^{\text{fam}}(\boldsymbol{q} , \boldsymbol{p}))$
parametrized by $U$
is locally constant.
\end{itemize}
Here $(l\times \mathrm{id})$ is a map from $U \times \mathbb{P}^1$ 
to $\mathrm{Sym}^3(\mathbb{C}^2) \times  (\mathbb{C}^*)^3 \times  \mathbb{P}^1$.
In general, this map \eqref{SolutionGeneral} is highly transcendental.
When the connections \eqref{OrigConn} have a special local formal data in \eqref{LocalFormalData},
we can construct the map \eqref{SolutionGeneral} algebraically.
For an algebraic isomonodromic family, 
$ \{ (q_j(\tilde{\boldsymbol{t}}) , \eta_j(\tilde{\boldsymbol{t}})) \}_{j=1,2,3}$ satisfies 
the irregular Garnier system \eqref{GarnierInSect2}.
So, if we have an algebraic isomonodromic family,
we obtain an algebraic solution of \eqref{GarnierInSect2}.

In particular,
we are interested in 
nonclassical algebraic solutions for irregular Garnier systems of rank $N>1$.
By \cite[Theorem 2]{DL1}, 
up to canonical transformations, 
there are exactly three nonclassical algebraic solutions for 
$N$-variable irregular Garnier systems with $N>1$.
The corresponding algebraic isomonodromic families are constructed by the pull-back method. 
In fact, the two of the three nonclassical algebraic solutions 
were already described in \cite[Section 9]{DL1} by this method. 
Remark that 
the irregular Garnier systems corresponding to these two solutions are $2$-variable.
So the explicit forms of the irregular Garnier systems are known (\cite{Kim} and \cite{Kawa}). 
On the other hand, for the last one of the three nonclassical algebraic solutions, 
we are able to compute the algebraic isomonodromic family by the pull-back method.
But the corresponding algebraic solution was not described in \cite{DL1},
since the explicit form of the irregular Garnier system had not been known in that case. 
But, now we have 
an explicit form of this irregular Garnier system, which is just \eqref{GarnierInSect2}.
So we will construct an algebraic isomonodromic family by using the pull-back method as in \cite{DL1}.
By this algebraic isomonodromic family, 
we will give 
the remaining nonclassical algebraic solution
of irregular Garnier systems.

\begin{Rem}
In \cite[Theorem 2]{DL1}, 
the three nonclassical algebraic solutions for 
$N$-variable irregular Garnier systems with $N>1$
are characterized by the following notation:
$$
\begin{pmatrix}
0&1&1 \\
\frac{1}{3} & 0 & 1
\end{pmatrix},\qquad 
\begin{pmatrix}
1&2\\
0 & 1
\end{pmatrix},\qquad 
\begin{pmatrix}
1&1&1 \\
0 & 0 & 1
\end{pmatrix}.
$$
Here, the explanation of the meaning of the notation is omitted
(see \cite[Section 1]{DL1}).
The irregular Garnier systems corresponding to the first one and second one 
are $2$-variable. 
These solutions were described in \cite{DL1}.
The third one is our algebraic solution in Theorem \ref{MainTheorem},
which is our main result.
\end{Rem}

\subsection{Isomonodromic family due to Diarra--Loray}

Now we construct the algebraic isomonodromic family 
corresponding to the remaining nonclassical algebraic solution
by the pull-back method.
The idea of the pull-back method has been used by Doran, Kitaev, Andreev, Vidunas, and Diarra 
(see \cite{Dor, AK, AK2, Kita, Kita1,Kita2, Kita3, VK,VK1, Dia}).
We will prepare a fixed system on $\mathbb{P}^1$ and an algebraic family of ramified covers
$\mathbb{P}^1\rightarrow \mathbb{P}^1$.
If we have them, then 
the pull-back of this fixed system by 
the ramified covers is an algebraic isomonodromic family
(for details, see the paragraphs after \cite[Theorem 1]{DL1}).

The fixed system on $\mathbb{P}^1$ and the algebraic family of ramified covers
which give the algebraic isomonodromic family corresponding to 
the remaining nonclassical algebraic solution are 
described in \cite{DL1} (see the first line of Table 2 in \cite{DL1}).
We recall them.
The fixed system on $\mathbb{P}^1$ is the following second-order scalar equation
\begin{equation}\label{FixedSystem}
\frac{d^2u}{dz^2} +\frac{2}{3\, z} \frac{du}{dz} -\frac{1}{z} u =0.
\end{equation}
This system has one logarithmic pole at $z=0$ and 
a ramified irregular singular point at $z=\infty$.
The algebraic family of ramified covers is 
$\phi_{\boldsymbol{s}} \colon \mathbb{P}^1 \rightarrow \mathbb{P}^1$ defined by
$$
z=\phi_{\boldsymbol{s}}(x) = 
{\frac { \left( s_{{2}}x \left( x-1 \right) +s_{{1}}x+s_{{0}} \left( 1
-x \right)  \right) ^{3}}{{x}^{2} \left( x-1 \right) ^{2}}}.
$$
This family $\phi_{\boldsymbol{s}}$ is parametrized by 
$\boldsymbol{s}=(s_0,s_1,s_2) \in (\mathbb{C}^*)^3$.

Now we consider the pull-back of the fixed system \eqref{FixedSystem} 
by the algebraic family $\phi_{\boldsymbol{s}}$ of ramified covers.
Then we have a second-order scalar equation.
We have a connection matrix in companion form corresponding to this scalar equation.
We may transform this connection into the following form:
\begin{equation}\label{IsomFami}
\begin{pmatrix}
0  & \frac{1}{x^2(x-1)^2} \\
\frac{ 4\, s_2 (s_2 x^2 -(s_0-s_1+s_2)x +s_0)Q(x; \boldsymbol{s})^2}{x^2(x-1)^2} &
-\frac{1}{3\, x}-\frac{1}{3(x-1)} - \sum_{i=1}^2 \frac{Q'(x; \boldsymbol{s})}{Q(x; \boldsymbol{s})}
\end{pmatrix} dx
\end{equation}
Here we set
$$
Q(x; \boldsymbol{s}) := x^3 +\frac{s_0-s_1-3\, s_2}{2\, s_2} x^2 
+\frac{-3\, s_0-s_1+ s_2}{2\, s_2}x +\frac{s_0}{s_2}.
$$
This connection matrix \eqref{IsomFami} is extended as
\begin{equation}\label{IsomFamiEx}
E_4 \longrightarrow E_4 \otimes \Omega^1_{\mathbb{P}^1} (D + [Q(x; \boldsymbol{s})=0]) 
\end{equation}
naturally.
Here $[Q(x; \boldsymbol{s})=0]$ is the effective divisor on $\mathbb{P}^1$ defined by
$Q(x; \boldsymbol{s})=0$.
The connection \eqref{IsomFamiEx} has poles at 
$0,1,\infty$ and at the zeros of the polynomial $Q(x; \boldsymbol{s})$ in $x$.
The poles at $0,1$, and $\infty$ are unramified irregular singular points of \eqref{IsomFami}.
If we define $\tilde{t}_{i}$ which satisfies ${\tilde{t}_{i}}^2  = s_i^3$
for each $i=0,1,2$, 
this connection 
has the local formal data as in \eqref{LocalFormalData}
with $t_0=\tilde{t}_0$, $t_1=\tilde{t}_1$, and $t_2=\tilde{t}_2$.
On the other hand, 
the zeros of $Q(x; \boldsymbol{s})$ are apparent singularities of \eqref{IsomFami}.
The corresponding parameters $p_j$ ($j=1,2,3$) 
appeared in \eqref{NormalForm} are zero: $p_1=p_2=p_3=0$.

We take an open set $U\subset (\mathbb{C}^*)^3$, which is small enough. 
We define a map
\begin{equation}\label{MapDL}
\begin{aligned}
l_{\text{DL}} \colon  U &\longrightarrow  
\mathcal{M} \times (\mathbb{C}^*)^3 \\
\tilde {\boldsymbol{t}}=(\tilde t_0,\tilde t_1,\tilde t_2) &\longmapsto
 ( \{ (q_j(\tilde{\boldsymbol{t}}) , \eta_j(\tilde{\boldsymbol{t}})) \}_{j=1,2,3},
 \tilde {\boldsymbol{t}})
\end{aligned}
\end{equation}
by taking branches of $\tilde{t}_{i}^2  = s_i^3$ $(i=0,1,2)$
and $Q(q_j; \boldsymbol{s}) =0$ $(j=1,2,3)$ for $\tilde{\boldsymbol{t}} \in U$.
Here $\eta_j$ ($j=1,2,3$) are determined by $\eta_j  = \frac{1}{3\, q_j} +\frac{1}{3(q_j-1)}$.

\begin{proof}[Proof of Theorem \ref{MainTheorem}]
We can check that $l_{\text{DL}}^* \nabla_4^{\text{fam}}(\boldsymbol{q}, \boldsymbol{p})$
coincides with \eqref{IsomFami}.
By the construction of \eqref{IsomFami},
the pull-back $l_{\text{DL}}^* \nabla_4^{\text{fam}}(\boldsymbol{q}, \boldsymbol{p})$
is isomonodromic. 
The local formal data is just \eqref{LocalFormalData}
with $t_0=\tilde{t}_0$, $t_1=\tilde{t}_1$, and $t_2=\tilde{t}_2$.
The functions $q_j ,\eta_j$ $(j=1,2,3)$ in \eqref{MapDL} 
are defined by 
\begin{equation*}
\left\{
\begin{aligned}
&q_j^3 +\frac{s_0-s_1-3\, s_2}{2\, s_2} q_j^2 
+\frac{-3\, s_0-s_1+ s_2}{2\, s_2}q_j +\frac{s_0}{s_2} =0 & j=1,2,3\\ 
&\eta_j=\frac{1}{3\, q_j} +\frac{1}{3\, (q_j-1)} & j=1,2,3 \\
&t_{i}^2  = s_i^3  & i=0,1,2.
\end{aligned}
\right.
\end{equation*}
implicitly.
So,
we have that $q_j ,\eta_j$ defined in \eqref{AlgebraicSol}
satisfy the irregular Garnier system defined in Definition \ref{Garnier_Def}.
\end{proof}

\subsection{The corresponding $\tau$-function}\label{2022_9_17_19_52}

We define functions $H_{t_i} (\boldsymbol{t})$ ($i=0,1,2$) by
$$
H_{t_i} (\boldsymbol{t}) = H_{t_i} (\boldsymbol{t} ,\boldsymbol{q}(\boldsymbol{t}), 
\boldsymbol{\eta}(\boldsymbol{t} ) ),
$$ 
where $\boldsymbol{q}(\boldsymbol{t}) =(q_1(\boldsymbol{t}),q_2(\boldsymbol{t}),q_3(\boldsymbol{t}))$ and 
$\boldsymbol{\eta}(\boldsymbol{t}) =(\eta_1(\boldsymbol{t}),\eta_2(\boldsymbol{t}),\eta_3(\boldsymbol{t}))$
are the algebraic solution
defined by \eqref{MainTheorem}, implicitly.
By the equalities \eqref{GarnierInSect2} and 
\eqref{2022_9_16_23_08},
we have that 
$$
\frac{\partial H_{t_i} (\boldsymbol{t})}{\partial t_j} =
\frac{\partial H_{t_j} (\boldsymbol{t} ,\boldsymbol{q}, \boldsymbol{\eta} ) }{
\partial t_i}|_{\boldsymbol{q} =\boldsymbol{q}(\boldsymbol{t}) ,
\boldsymbol{\eta} =\boldsymbol{\eta}(\boldsymbol{t}) }.
$$
On the other hand, 
by computation using 
the explicit formula \eqref{Hamiltonians} and
using some computing environment (for example {\it Maple}),
we have that
\begin{equation*}\label{2022_9_15_23_16}
\{ H_{t_i} (\boldsymbol{t} ,\boldsymbol{q}, \boldsymbol{\eta} ), 
H_{t_j} (\boldsymbol{t} ,\boldsymbol{q}, \boldsymbol{\eta} )  \}|_{\boldsymbol{q} 
=\boldsymbol{q}(\boldsymbol{t}) ,
\boldsymbol{\eta} =\boldsymbol{\eta}(\boldsymbol{t}) }=0.
\end{equation*}
By this equality and the equality \eqref{2022_9_16_23_08}, we have that 
$$
\begin{aligned}
\frac{\partial H_{t_j} (\boldsymbol{t})}{\partial t_i} 
&=
\frac{\partial H_{t_i} (\boldsymbol{t} ,\boldsymbol{q}, \boldsymbol{\eta} ) }{
\partial t_j}|_{\boldsymbol{q} =\boldsymbol{q}(\boldsymbol{t}) ,
\boldsymbol{\eta} =\boldsymbol{\eta}(\boldsymbol{t}) } \\
&=
\frac{\partial H_{t_j} (\boldsymbol{t} ,\boldsymbol{q}, \boldsymbol{\eta} ) }{
\partial t_i}|_{\boldsymbol{q} =\boldsymbol{q}(\boldsymbol{t}) ,
\boldsymbol{\eta} =\boldsymbol{\eta}(\boldsymbol{t}) }
=\frac{\partial H_{t_i} (\boldsymbol{t})}{\partial t_j}.
\end{aligned}
$$
Then if 
we define a 1-form $\varpi$ as
$$
\varpi =H_{t_0} (\boldsymbol{t}) dt_0 + H_{t_1} (\boldsymbol{t}) dt_1+ H_{t_2} (\boldsymbol{t}) dt_2,
$$
then this 1-form is closed.
Since $\varpi $ is closed, we can define the {\it $\tau$-function}
as in \cite[Section 5]{JMU}.
That is, the function $\tau$ is defined by the equation $d \ln \tau =\varpi$.

Now we will give explicit description of the $\tau$-function.
If we compute $\varpi$ by using the explicit formula \eqref{Hamiltonians} and
using some computing environment, then we have the following explicit description of $\varpi$:
$$
\begin{aligned}
\varpi&={\frac {1+36\,{s_{{0}}}^{3}-216\, \left( s_{{1}}+\,s_{{2}}
 \right) {s_{{0}}}^{2}-108\, \left( s_{{1}}-s_{{2}} \right) ^{2}s_{{0}
}}{24\,s_{{0}}}}ds_0\\
&\quad +{\frac {1+36\,{s_{{1}}}^{3}-216\, \left( s_{{0}}+ s_{{2}}
 \right) {s_{{1}}}^{2}-108\, \left( s_{{0}}-s_{{2}} \right) ^{2}s_{{1}
}}{24\,s_{{1}}}}ds_1\\
&\quad +{\frac {1+36\,{s_{{2}}}^{3}-216\, \left( s_{{0}}+ s_{{1}}
 \right) {s_{{2}}}^{2}-108\, \left( s_{{0}}-s_{{1}} \right) ^{2}s_{{2}
}}{24\,s_{{2}}}}ds_2 .
\end{aligned}
$$
By this formula, we have that the 1-form $\varpi$ is exact.
Indeed, if we set 
$$
\begin{aligned}
F(s_0,s_1,s_2) &:= 
{\frac {{s_{{0}}}^{3} + {s_{{1}}}^{3}+ {s_{{2}}}^{3}}{2}}
+{\frac {\ln  \left( s_{{0}}s_{{1}}s_{{2}} \right) }{24}} \\
&\quad -{\frac {9 \, \left( \left( s_{{1}}-s_{{2}} \right) ^{2}s_{{0}} 
+\left( s_{{0}}-s_{{2}} \right) ^{2}s_{{1}} 
+ \left( s_{{1}}-s_{{0}} \right) ^{2}s_{{2}}   \right) }{2} } 
-18\,s_{{0}}s_{{1}}s_{{2}},
\end{aligned}
$$
then $\varpi = dF(s_0,s_1,s_2)$.
So the $\tau$-function has the following description 
$$
 \tau = c \cdot e^{F(s_0,s_1,s_2)},
$$
where $c$ is a constant.

\section*{Appendix}

Our irregular Garnier system is described explicitly 
and these functions $q_j ,\eta_j$ have simple description.
So we may recheck 
these functions $q_j ,\eta_j$ satisfy the irregular Garnier system
by direct calculation.
Now we will recheck 
that $q_j ,\eta_j$ ($j=1,2,3$) which are implicitly defined by \eqref{AlgebraicSol}
satisfy the equation
$\frac{\partial q_j}{\partial t_i} =  \frac{\partial H_{t_i}(\boldsymbol{t},\boldsymbol{q},\boldsymbol{\eta})}{\partial \eta_j}$
by direct calculation
using some computing environment (for example {\it Maple}).
We may also recheck 
that $q_j ,\eta_j$ satisfy the equation
$\frac{\partial \eta_j}{\partial t_i} = -  \frac{\partial H_{t_i}(\boldsymbol{t},\boldsymbol{q},\boldsymbol{\eta})}{\partial q_j}$
by direct calculation.
But the calculation is too huge to describe here.
So we omit the calculation of 
$\frac{\partial \eta_j}{\partial t_i} = -  \frac{\partial H_{t_i}(\boldsymbol{t},\boldsymbol{q},\boldsymbol{\eta})}{\partial q_j}$.

We consider the following composition
$$
\begin{aligned}
\mathcal{M} \subset 
\mathrm{Sym}^3(\mathbb{C}^2) &\longrightarrow
\mathrm{Sym}^3(\mathbb{C}) \longrightarrow  \mathbb{C}^3 \\
\{ (q_1,\eta_1), (q_2,\eta_2), (q_3,\eta_3)\} &\longmapsto
\{ q_1, q_2,q_3 \} \longmapsto  (u_1(\boldsymbol{q}),u_2(\boldsymbol{q}),u_3(\boldsymbol{q})),
\end{aligned}
$$
where we set 
\begin{equation*}
\left\{
\begin{aligned}
&\sigma_1(\boldsymbol{q}):=q_1+q_2+q_3 \\
&\sigma_2(\boldsymbol{q}):=q_1q_2 +q_1q_3 +q_2q_3  \\
&\sigma_3(\boldsymbol{q}):=q_1q_2q_3 .
\end{aligned}
\right.
\end{equation*}
Now we consider a transformation of the equation
$\frac{\partial q_j}{\partial t_i} 
=  \frac{\partial H_{t_i}(\boldsymbol{t},\boldsymbol{q},\boldsymbol{\eta})}{\partial \eta_j}$
via this composition.
We will give a system on $\mathbb{C}^3$.
First, we define functions on $\mathrm{Sym}^3(\mathbb{C}^2)$ as follows:
$$
\begin{aligned} 
K^{(1)}_{t_i} (\boldsymbol{q},\boldsymbol{\eta})
&:= \frac{\partial H_{t_i}}{\partial \eta_1}+\frac{\partial H_{t_i}}{\partial \eta_2}
+\frac{\partial H_{t_i}}{\partial \eta_3} \\
K^{(2)}_{t_i} (\boldsymbol{q},\boldsymbol{\eta})
&:= (q_2+q_3)\frac{\partial H_{t_i}}{\partial \eta_1}
+(q_1+q_3)\frac{\partial H_{t_i}}{\partial \eta_2}
+(q_1+q_2)\frac{\partial H_{t_i}}{\partial \eta_3} \\
K^{(3)}_{t_i} (\boldsymbol{q},\boldsymbol{\eta})
&:=q_2q_3\frac{\partial H_{t_i}}{\partial \eta_1}
+q_1q_3\frac{\partial H_{t_i}}{\partial \eta_2}
+q_1q_2\frac{\partial H_{t_i}}{\partial \eta_3}.
\end{aligned}
$$
Second, we define a section $\sigma$ of the projection 
$\mathrm{Sym}^3(\mathbb{C}^2) \rightarrow\mathrm{Sym}^3(\mathbb{C})$ as
$$
\begin{aligned}
\sigma \colon \mathrm{Sym}^3(\mathbb{C}) &\longrightarrow
\mathrm{Sym}^3(\mathbb{C}^2)  \\
\{ q_1, q_2,q_3 \}  &\longmapsto
\{ (q_1,\eta_1(\boldsymbol{q})), (q_2,\eta_2(\boldsymbol{q})), (q_3,\eta_3(\boldsymbol{q}))\}
\end{aligned}
$$
Here $\eta_j(\boldsymbol{q})$ is defined by
$\eta_j(\boldsymbol{q})=\frac{1}{3q_j} +\frac{1}{3(q_j-1)} $
for each $j=1,2,3$.
Remark that this formula appeared in \eqref{AlgebraicSol}.
Third, we consider the pull-back of $K^{(j)}_{t_i}$ by $\sigma$
for each $j=1,2,3$ and $i=0,1,2$.
Then we have the following transformation of the equation
$\frac{\partial q_j}{\partial t_i} =  \frac{\partial H_{t_i}(\boldsymbol{t},\boldsymbol{q},\boldsymbol{\eta})}{\partial \eta_j}$:
$$
\frac{\partial \sigma_j}{\partial t_i} =
\sigma^* K^{(j)}_{t_i} (\boldsymbol{q},\boldsymbol{\eta})
$$
for any $i=0,1,2$ and $j=1,2,3$.
We can calculate $\sigma^* K^{(j)}_{t_i} (\boldsymbol{q},\boldsymbol{\eta})$ explicitly
for any $i=0,1,2$ and $j=1,2,3$.
These pull-backs induce functions on $\mathbb{C}^3$.
Finally, we have the following system: 
\begin{equation}\label{reductionQ}
\left\{
\begin{aligned}
&\frac{\partial \sigma_1}{\partial t_0} =  {\frac {\sigma_3}{3\,t_{{0}}}} \\
&\frac{\partial \sigma_2}{\partial t_0} ={\frac {\sigma_3 }{\,t_{{0}}}} \\
&\frac{\partial \sigma_3}{\partial t_0} ={\frac {2\, \sigma_3 }{3\,t_{{0}}}}
\end{aligned}  \right. \qquad 
\left\{
\begin{aligned}
&\frac{\partial \sigma_1}{\partial t_1} ={\frac { \sigma_3 -\sigma_2 +\sigma_1 -1}{3\,t_{{1}}}} \\ 
&\frac{\partial \sigma_2}{\partial t_1} =-{\frac { \sigma_3 -\sigma_2 +\sigma_1 -1}{3\,t_{{1}}}} \\
&\frac{\partial \sigma_3}{\partial t_1} =0
\end{aligned} \right. \qquad 
\left\{
\begin{aligned}
&\frac{\partial \sigma_1}{\partial t_2} =-{\frac {2\,\sigma_1-3}{3\,t_{{2}}}} \\
&\frac{\partial \sigma_2}{\partial t_2} =-{\frac {2\,  \sigma_2 -1}{3\,t_{{2}}}} \\
&\frac{\partial \sigma_3}{\partial t_2} =-{\frac {2\,  \sigma_3 }{3\,t_{{2}}}}
\end{aligned}\right.
\end{equation}

We can check directly that 
$\sigma_1= - \frac{s_0 - s_1 - 3 s_2}{2s_2}$,
$\sigma_2= \frac{-3s_0 - s_1+s_2}{2s_2}$, 
$\sigma_3=- \frac{s_0}{s_2}$,
and $t_i^2 =s_i^3$ ($i=0,1,2$)
satisfy the system \eqref{reductionQ}.
By the construction of the section $\sigma$, 
we have that $q_j ,\eta_j$ ($j=1,2,3$) which are implicitly defined by \eqref{AlgebraicSol}
satisfy the equation
$\frac{\partial q_j}{\partial t_i} =  \frac{\partial H_{t_i}(\boldsymbol{t},\boldsymbol{q},\boldsymbol{\eta})}{\partial \eta_j}$.

\subsection*{Acknowledgments}

The author would like to thank Professor Frank Loray 
for leading him to the subject treated in this paper and also for valuable discussions.
He would like to thank Professor Ryo Ohkawa 
and Takafumi Matsumoto for valuable discussions.
He is supported by Japan Society for the Promotion of Science KAKENHI 
Grant Numbers 22H00094 and 19K14506.
He is very grateful to the anonymous referee's insightful 
suggestions which helped to improve the paper.
In particular, Section \ref{2022_9_17_19_52}
has been added 
in response to referee's comment.

\noindent
Komyo Arata\\
Center for Mathematical and Data Sciences\\
Kobe University\\
1-1 Rokkodai-cho, Nada-ku, Kobe, 657-8501\\
Japan\\
akomyo@math.kobe-u.ac.jp

\end{document}